\documentclass[journal,transmag]{IEEEtran}
\usepackage[utf8]{inputenc}
\usepackage[top=0.7in, left=0.65in]{geometry}
\usepackage{color}
\usepackage{graphicx}
\usepackage{bm}
\usepackage{amsmath, amssymb}
\usepackage{siunitx}
\usepackage{placeins}
\usepackage{todonotes}

\usepackage{tikz}
\usepackage{pgfplots}
\pgfplotsset{compat=newest}
\usepgfplotslibrary{groupplots}
\usepackage{tikzscale}
\definecolor{color0}{RGB}{237, 28, 36}
\definecolor{color1}{HTML}{50B695}
\definecolor{color2}{HTML}{EE7A34}
\definecolor{color3}{HTML}{0083CC}
\definecolor{colorx}{HTML}{951169}

\usepackage[T1]{fontenc}
\graphicspath{{graph/}{data/}}

\begin{document}
\title{\fontsize{17pt}{17pt}\selectfont A Parallel-In-Time Adjoint Sensitivity Analysis for a B6 Bridge-Motor Supply Circuit}
\author{\IEEEauthorblockN{Julian Sarpe\IEEEauthorrefmark{1},
Andreas Klaedtke\IEEEauthorrefmark{2}, and
Herbert De Gersem\IEEEauthorrefmark{3}}
\vspace{5pt}
\IEEEauthorblockA{\IEEEauthorrefmark{1}
TU Darmstadt,
Schlo{\ss}gartenstr.~8, 64289 Darmstadt, 
Germany, julian\_johannes.buschbaum@tu-darmstadt.de}
\IEEEauthorblockA{\IEEEauthorrefmark{2}
Robert Bosch GmbH,
Robert-Bosch-Campus 1, 71272 Renningen, Germany, andreas.klaedtke@de.bosch.com}
\IEEEauthorblockA{\IEEEauthorrefmark{3}
TU Darmstadt,
Schlo{\ss}gartenstr.~8, 64289 Darmstadt,  Germany, herbert.degersem@tu-darmstadt.de}
}

\IEEEtitleabstractindextext{%
\begin{abstract}
	This paper presents a parallel-in-time adjoint sensitivity analysis which combines a transient adjoint sensitivity analysis with the parareal approach in order to significantly speed up the simulation.
	The adjoint method is the method of choice to calculate the sensitivities in a many-parameter setting. 
	In order to obtain sensitivity information that is time-dependent, multiple adjoint problems must be solved.
	This slows down the simulation wall-clock time and leaves a large optimization potential for the analysis.
	The parareal is applied to the adjoint solution, significantly speeding up the adjoint solution for every timestep respectively.
\end{abstract}

\begin{IEEEkeywords}
Adjoint method, Modified nodal analysis, Parareal, Power-electronic devices
\end{IEEEkeywords}}

\maketitle
\thispagestyle{empty}
\pagestyle{empty}

\section{Introduction}
	Sensitivity analysis is a vital task during the optimization stage of the design process. 
	It is also a fundamental step in finding the root cause of issues within an existing design.
	Having the sensitivities of multiple design parameters is beneficial in a statistical tolerance analysis \cite{article:adjointtd}, \cite{article:cao2003adjoint}.
	Here, we are concentrating on methods for the sensitivity analysis w.r.t. all parameters for large scale time-periodic circuits,  
	which can be efficiently obtained along the adjoint formulation of the problem.
	However, time-dependent sensitivities require the solution of numerous adjoint problems, which cause a computational bottleneck. 
	Parallel algorithms cut down this wall-clock time. 
	This paper proposes to combine the adjoint solver computation with the parareal algorithm~\cite{article:parareal} to improve the wall clock time.
	While split time domain methods are well known in literature~\cite{article:checkpointAdjoint}, parallel-in-time methods have not been applied to adjoint sensitivity analysis to the authors' knowledge.

\section{Adjoint Sensitivity Analysis in Time Domain}
	The modified nodal analysis (MNA) of a circuit in time domain leads to a system of differential algebraic equations (DAE)~\cite{article:adjointtdder} of the form
	\begin{equation}
		\bm{F}(\bm{\varphi},\dot{\bm{\varphi}},t) = \bm{J}_\mathrm{C} \dot{\bm{\varphi}} + \bm{J}_\mathrm{G} \bm{\varphi} - \bm{i}_\mathrm{s} = 0,
		\label{eq:DAE}
	\end{equation}
	where $\bm{J}_\mathrm{C} = \bm{Y}_\mathrm{C} - \partial\bm{i}_\mathrm{nl}/\partial \dot{\bm{\varphi}}$ contains capacitor and inductor stamps, $\bm{J}_\mathrm{G} = \bm{Y}_\mathrm{G} - \partial\bm{i}_\mathrm{nl}/\partial \bm{\varphi}$ source and conductance stamps, $\bm{i}_\mathrm{s}$ represents the independent sources and $\bm{i}_\mathrm{nl}$ the nonlinear contributions. $\varphi(t)$ gathers the degrees of freedom (DoFs) (nodal voltages as well as currents) from which the quantities of interest $U(\varphi, t)$ are derived.
	\par
	In order to obtain the sensitivities, Eq.~\eqref{eq:DAE} is derived w.r.t. the parameters $p_i$. The parameters are related to components in the network, such as RLC-elements. 
	Introducing the adjoint variable $\lambda$ and integrating over the time interval [0, $t_m$]~\cite{article:cao2003adjoint} gives the equation
	\begin{multline}
		\int_{0}^{t_m}\lambda^\mathsf{T}\left(\bm{J}_\textrm{C} \frac{\mathrm{d} \dot{\bm{\varphi}}}{\mathrm{d} p_i}(t) + \bm{J}_\textrm{G} \frac{\mathrm{d} \bm{\varphi}}{\mathrm{d} p_i}(t)\right)\mathrm{d}t \\
		+ \int_{0}^{t_m}\lambda^\mathsf{T}\left(\frac{\mathrm{d} \bm{J}_\textrm{C}}{\mathrm{d} p_i}\dot{\bm{\varphi}}(t) + \frac{\mathrm{d} \bm{J}_\textrm{G}}{\mathrm{d} p_i}\bm{\varphi}(t)\right)\mathrm{d}t = 0.
		\label{eq:theprob}
	\end{multline}
	The unknown terms $\mathrm{d}\dot{\bm{\varphi}}/\mathrm{d}p$ and $\mathrm{d}\bm{\varphi}/\mathrm{d}p$ are eliminated by appropriate choice of $\lambda$
	turning out to be the solution of the time reversed DAE with the initial condition $\lambda(t=t_m) = 0$.

	\begin{equation}
		\bm{J}_\textrm{C}^\mathsf{T} \dot{\lambda} - \bm{J}_\textrm{G}^\mathsf{T} \lambda = \frac{\partial U}{\partial \bm{\varphi}}.
	\label{eq:adjoint}
	\end{equation}
	Using the integration-by-parts approach for the first term in~\eqref{eq:adjoint} similar to a variational formulation~\cite{article:adjointtdder} results in
	\begin{equation}
		\int_0^{t_m} \frac{\mathrm{d} U}{\mathrm{d} p_i} ~\mathrm{d} t = \int_0^{t_m} \lambda^\mathsf{T} \left(\frac{\mathrm{d} \bm{J}_\textrm{C}}{\mathrm{d} p_i}\dot{\bm{\varphi}}(t) + \frac{\mathrm{d} \bm{J}_\textrm{G}}{\mathrm{d} p_i}\bm{\varphi}(t)\right)\mathrm{d} t.
		\label{eq:sensintegral}
	\end{equation}
	To obtain the sensitivity for a quantity of interest (QoI) at a certain time instant, the generalized version of Leibniz integral rule can be applied~\cite{article:cao2003adjoint} and the sensitivity at time $t_m$ follows as
	\begin{multline}
		\frac{\mathrm{d} U}{\mathrm{d} p_i} (t_m) = \\ \int_0^{t_m} \left(\frac{\partial}{\partial t_m} \lambda^\mathsf{T}\right) \left(\frac{\mathrm{d} \bm{J}_\textrm{C}}{\mathrm{d} p_i}\dot{\bm{\varphi}}(t) + \frac{\mathrm{d} \bm{J}_\textrm{G}}{\mathrm{d} p_i}\bm{\varphi}(t)\right)\mathrm{d}t.
		\label{eq:sensSolIntegral}
	\end{multline}
	Note here that individual adjoint solutions are required for each analyzed time instant, because of the boundary condition required to solve~\eqref{eq:adjoint} and~\eqref{eq:sensSolIntegral} as described in~\cite{article:cao2003adjoint}. 
	As a result, the adjoint system needs to be solved individually for each time instant analyzed.
	This requires a large number of numerical simulations, which will be optimized and improved with the method proposed in this paper.

\setlength{\tabcolsep}{4pt} 
\section{Parareal}
\begin{figure}[t]
		\centering
		\includegraphics[width=.65\linewidth]{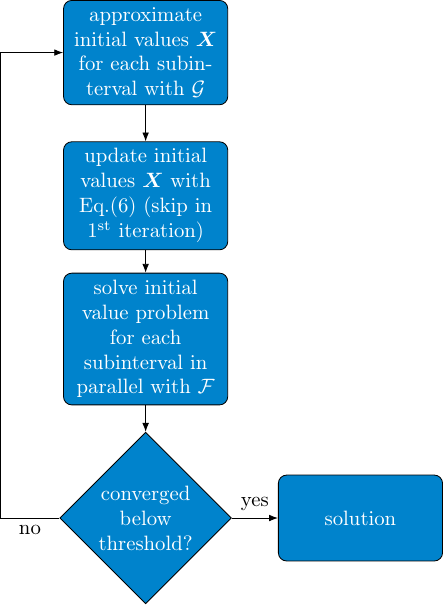}
		\caption{Flowchart outlining the iteration steps along the parareal approach.}
		\label{fig:flowchart}
	\end{figure}

	Parareal~\cite{article:parareal} approaches are based on the idea that DAEs can be solved by splitting an interval into $N$ subintervals and then solving these $N$ initial value problems in parallel. 
	A fast coarse solver  $\mathcal{G}$ delivers the initial values $\bm{X}$ for $N$ temporal subintervals in each iteration.
	The fine solvers $\mathcal{F}$, running in parallel, provide updates for the Newton iterations converging to the
	solution. The process is outlined as a flowchart in Fig.~\ref{fig:flowchart}.
	Both fine solver $\mathcal{F}$ and coarse solver $\mathcal{G}$ are black-boxes within the simulation and can be chosen accordingly. 
	The coarse solver is usually a fast but less accurate solver, such as a transient solver with large time steps. 
	The primary function of the coarse solver is to offer approximate initial values for each subinterval. 
	Since it directly affects the serial part of the parallel algorithm, it is crucial for the coarse solver to be as efficient as possible (Amdahl's law~\cite{article:gustafson1988}).
	The fine solver $\mathcal{F}$ is an accurate but more expensive solver and is used to solve the initial value problem for each of the $N$ subintervals in parallel.
	The first iteration starts with approximate initial values. 
	This leads to discontinuities at the interfaces, which are minimized iteratively~\cite{article:pararealInd2}.
	The initial values for the subintervals are updated at the end of each iteration by:
	\begin{equation}
		\bm{X}^{k+1}_{n} = \mathcal{F}(\bm{X}^k_{n-1}) + \mathcal{G}(\bm{X}^{k+1}_{n-1}) - \mathcal{G}(\bm{X}^{k}_{n-1}).
		\label{eq:update}
	\end{equation}
	Here, 
	$n$ is the index of the subinterval starting point, 
	$k$ is the index of the iteration
	 and $\bm{X}$ are the initial values at the respective interfaces. 
	The iterative process ends when the discontinuities fall below a given error threshold.
	\par
	In the context of the adjoint sensitivity analysis, parareal can be applied to the forward transient circuit simulation as well as to the adjoint simulation. 
	
\section{Applications}
\subsection{Half-Wave Rectifier}
\begin{figure}[t]
	\centering
	\includegraphics[width=.6\linewidth]{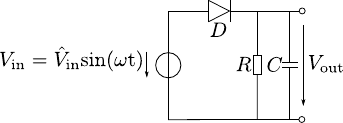}
	\caption{Half-wave rectifier functional circuit with sinusoidal input voltage $V_\mathrm{in}$.}
	\label{fig:HWrectcirc}
\end{figure}

\begin{figure}[t]
	\centering
	\includegraphics[width=.85\linewidth]{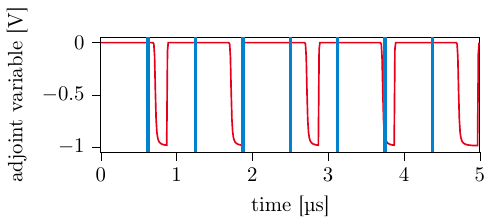}
	\caption{Adjoint variable $\lambda$ for the half-wave rectifier model in its first five periods of operation. The blue lines indicate the subinterval boundaries.}
	\label{fig:adjointvar}
\end{figure}

\begin{figure}[t]
	\centering
	\includegraphics[width=.9\linewidth]{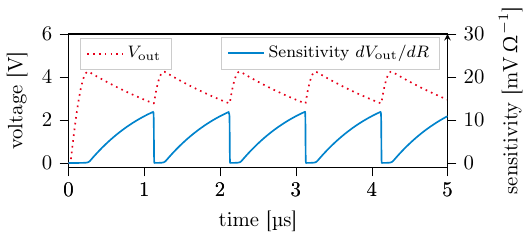}
	\caption{Sensitivity $\mathrm{d}V_\mathrm{out} / \mathrm{d}R$ for the half-wave rectifier model in its first five periods of operation.}
	\label{fig:sensHWRECT}
\end{figure}

A half-wave rectifier consists of only a diode, a resistor and a capacitor and a voltage source (Fig.~\ref{fig:HWrectcirc}). 
The adjoint variable $\lambda$ is calculated for each analyzed time instant as an initial value problem.
The adjoint variable $\lambda$ with initial value 0 at time instant $t_m$ calculated with the parareal method is shown in Fig.~\ref{fig:adjointvar}. 
The transient adjoint solutions are then plugged into Eq.~\ref{eq:sensSolIntegral} in order to obtain the sensitivity for each time instant $t_m$.
The resulting sensitivity is a function of time (as is the QoI). 
Fig.~\ref{fig:sensHWRECT} shows the sensitivity for the output voltage $V_\mathrm{out}$ w.r.t. the resistor $R$.
As expected, the resistor has high impact on the output voltage $V_\mathrm{out}$ while the capacitor is discharged, i.e. when the slope of $V_\mathrm{out}$ is negative, and is less influential otherwise.

\subsection{B6 Bridge-Motor Supply}
\begin{figure}[b]
	\centering
	\includegraphics[width=\linewidth]{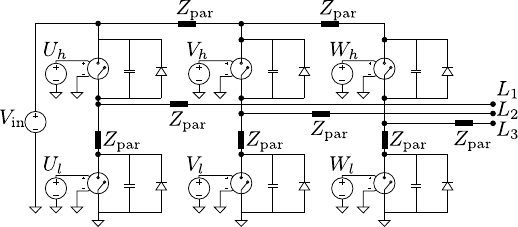}
	\caption{B6 bridge-motor supply circuit including parasitic EEC impedances $Z_\mathrm{par}$. $U_h$, $V_h$, $W_h$, $U_l$, $V_l$ and $W_l$ are PWM input voltages respectively.}
	\label{fig:B6circuit}
\end{figure}

To highlight the advantage of the adjoint approach, we study the sensitivities of a B6 bridge-motor supply (Fig.~\ref{fig:B6circuit}). 
The functional circuit is extended by an equivalent electrical circuit (EEC) of parasitic elements $Z_\mathrm{par}$, which arise from the field coupling between neighboring functional components and interconnects and consist of several hundred lumped circuit elements each. 
Each of the $Z_\mathrm{par}$ elements consists of a subcircuit with multiple EEC elements which are also coupled with each other, leading to a large-scale circuit problem. 
The switches are driven by the pulse wave modulated (PWM) signals $U$, $V$ and $W$ which each consist of a high-side and an opposing low-side.
The circuit behavior is also broadband due to the properties of the B6-bridge in general.
\par
\begin{figure}[t]
	\centering
	\includegraphics[width=.8\linewidth]{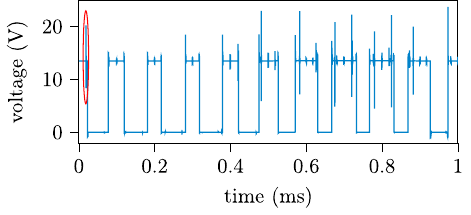}
	\caption{The QoI plotted over time for the first PWM-cylce period of the switch voltage.}
	\label{fig:solfull}
\end{figure}
\begin{figure}[b]
	\centering
	\includegraphics[width=.8\linewidth]{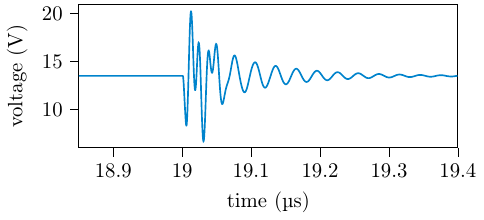}
	\caption{The QoI plotted over time for the first switching overshoot.}
	\label{fig:solzoomed}
\end{figure}

The QoI is chosen as the $L_1$-phase high-side switch voltage. 
This switch voltage (Fig.~\ref{fig:solfull}) contains the entire PWM sequence as well as high frequency overshoot oscillations. 
Thus, the switch voltage shows both the broadband and steep transient properties of the circuit behavior. 
The sensitivity is analyzed for the first switching overshoot, as indicated in Fig.~\ref{fig:solfull}.
Zooming in shows the shape of the relaxation oscillation, depicted in Fig.~\ref{fig:solzoomed}.
This switching induced relaxation oscillation is triggered by the switching event of the $L_2$-phase.
As a result, the mean voltage is nearly constant at $\SI{12.5}{\volt}$ in the observed period.
\par
\begin{figure}[h]
	\centering
	\includegraphics[width=.85\linewidth]{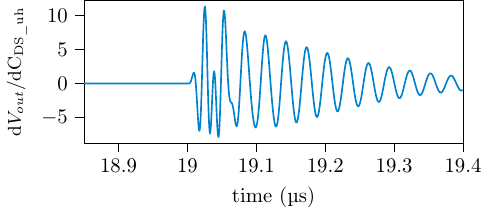}
	\caption{Sensitivity of output voltage w.r.t. switch capacitance of the upper L3-Phase $\mathrm{d}V_{out}/\mathrm{d}\mathrm{C_{DS\_uh}}$.}
	\label{fig:sensTDCDS}
\end{figure}
\begin{figure}[h]
	\centering
	\includegraphics[width=.95\linewidth]{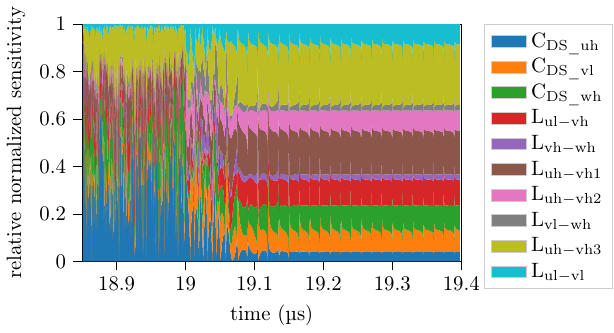}
	\caption{Relative normalized sensitivity for the ten most influential system parameters plotted over time.}
	\label{fig:sensTD}
\end{figure}

The sensitivity in the interval from $\SI{18.85}{\micro\second}$ to $\SI{19.4}{\micro\second}$ is analyzed w.r.t. all system parameters for all time instants, here chosen identically to the timesteps of the forward solver.
The sensitivity w.r.t. the switch capacitances directly follows the shape of the QoI. This behavior is depicted in Fig.~\ref{fig:sensTDCDS}, exemplarily for the $L_1$-phase high-side capacitance $\mathrm{C_{DS\_uh}}$.
The ten most influential system parameters are shown in Fig.~\ref{fig:sensTD}, represented by their relative sensitivity normalized w.r.t. the total sensitivity of the most influential parameters.
\par
Connecting the results to their physical meaning is important for interpretation.
Here, this interpretation will be demonstrated for two elements: $\mathrm{L_{uh-vh3}}$ and $\mathrm{C_{DS\_uh}}$.
The inductivity $\mathrm{L_{uh-vh3}}$ is associated to the interconnection between the $L_1$-phase high side switch $U_h$ and $L_2$-phase high side switch $V_h$. 
Therefore, $\mathrm{L_{uh-vh3}}$ can serve as a bypass path to the switch voltage. 
Furthermore, $\mathrm{L_{uh-vh3}}$ highly influences the oscillation behavior because it interacts with the switch capacitances of the $L_1$- and $L_2$-phases. 
The capacitance $\mathrm{C_{DS\_uh}}$ is particularly influential before the first switching. 
Afterwards, the relative influence is lower as the oscillation over the $L_1$ high side switch is mainly determined by the $L_2$-phase.
This is explained by the fact that the observed overshoot is triggered by the switching event of the $L_2$-phase.
\par
In many cases, the interpretation of results is easier to perform when they are available in frequency domain.
However, for non-periodic results, standard FFT is not suited to determine the power spectrum for the time series.
For non-periodic time series, the power spectrum can be estimated by a periodogram or similar methods.
The power spectrum for the sensitivities is estimated for the interval using Welchs method~\cite{article:welch}.
\begin{figure}[b]
	\centering
	\includegraphics[width=.95\linewidth]{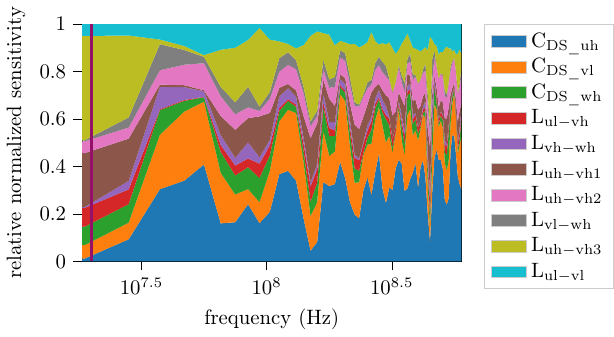}
	\caption{Power spectrum of the relative normalized sensitivity for the ten most influential system parameters. The vertical line ({\color{colorx} \textbf{\textemdash}}) indicates the relaxation oscillation frequency.}
	\label{fig:sensAC}
\end{figure}
The frequency domain graph shows the previously mentioned properties more clearly.
The inductivity $\mathrm{L_{uh-vh3}}$ strongly influences the harmonic, lower frequency oscillation particularly after $\SI{19.1}{\micro\second}$. 
The influence of $\mathrm{L_{uh-vh3}}$ on higher frequency oscillations is comparably weaker.
The capacitance $\mathrm{C_{DS\_uh}}$ weakly influences the lower frequency oscillation but has significant influence on high frequency components of the signal.
These results are particularly valuable in the optimization stage to improve the circuit behavior and prevent possible malfunctions at specific frequencies.

\section{Efficiency of the Parallel-in-Time Simulation}
The transient adjoint sensitivity analysis benefits from short time intervals of interest with as few analyzed time instances as possible.
The complexity of the simulation increases as the number of time instances to be analyzed grows, as the adjoint solution must be individually solved for every time instance by solving Eq.~\eqref{eq:adjoint}.
However, the simulation is sped up by parallelization. 
Here, the parareal approach is used to parallelize the solution of each individual adjoint problem to be solved at every time step.
The total wall-clock time for the simulation of one adjoint solution is measured to assess the efficiency of the parareal method.
\par
\begin{figure}[tb]
	\centering
	\includegraphics[width=.85\linewidth]{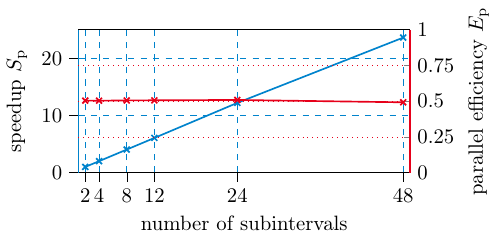}
	\caption{Speedup and parallel efficiency for different numbers of subintervals.}
	\label{fig:perform}
\end{figure}

The solver utilized for the fine solution is identical to the sequential solver.
The coarse solver is similar, but it solves only every 100th timestep, giving a very fast but less accurate result.
The parareal algorithm was observed to converge in two Newton iterations for this setup, regardless of the number of subintervals.
For the benchmark, the adjoint solution for $\SI{19.1}{\micro\second}$ is used as reference, since it is in the middle of the observed time period. Table~\ref{tab:sim_duration} shows the simulation duration for different numbers of subintervals.
\begin{table}[tb]
	\centering
	\caption{Simulation duration for different numbers of subintervals.}
	\begin{tabular}{|l|l|l|l|}
		\hline
		subintervals \hspace{-.6em} & fine solution (s) \hspace{-.6em} & coarse solution (s) \hspace{-.6em} & total wall-clock time (s) \hspace{-.6em} \\ \hline\hline
		2 & 25.07 & 0.24 & 50.62 \\ \hline
		4 & 12.55 & 0.12 & 25.34 \\ \hline
		8 & 6.25 & 0.06 & 12.62 \\ \hline
		12 & 4.16 & 0.04 & 8.40 \\ \hline
		24 & 2.07 & 0.02 & 4.18 \\ \hline
		48 & 1.07 & 0.01 & 2.16 \\ \hline
		\multicolumn{3}{|l|}{sequential solution: } & 51.08 \\ \hline
	\end{tabular}
	\label{tab:sim_duration}
\end{table}
The total parallel wall-clock time $T_{\mathrm{p}}$ is compared to total wall-clock time $T_\mathrm{s}$ of the sequential solver without parareal by the speedup
\begin{equation}
	S_{\mathrm{p}} = T_\mathrm{s} / T_{\mathrm{p}}.
 	\label{eq:speedup}
\end{equation}
This comparison shows the improvement of the parallel-in-time adjoint sensitivity analysis and highlights its efficiency.
The calculation of parallel efficiency 
\begin{equation}
E_{\mathrm{p}} = S_{\mathrm{p}} / N
\label{eq:par_eff}
\end{equation}
allows for quantification of the impact of parallelization on the overall efficiency of the adjoint sensitivity analysis.
Both parallel efficiency and speedup are shown in Fig.~\ref{fig:perform}.
The speedup grows linearly as the number of subintervals increases.
Consequently, the parallel efficiency remains constant regardless of the number of subintervals.
The linear speedup is explained by the negligible share of the coarse solver on the overall runtime and also because no communication between the parallel instances of parareal are necessary.
As a result, the algorithm scales very well with additional number of computational units (threads, tasks, processors).
Parareal is therefore a well suited measure to speed up the adjoint sensitivity analysis in the given application case.
\FloatBarrier

\section{Conclusions and Outlook}
Parallelization in time is a well suited approach to optimize the performance of adjoint sensitivity analysis. 
The simulation can be sped up linearly with the number of subintervals.
The parallel efficiency is constant for all considered numbers of subintervals. 
Consequently, the parareal approach can be used to greatly decrease simulation times, particularly when a massively parallel computing system is available, such as with simulation clusters.
Further improvements can be made by applying the parareal both to the forward transient solution as well as the adjoint solution.

\end{document}